\documentclass[12pt]{article}
\usepackage{setspace}
\usepackage[hmargin=3 cm,vmargin={2.5 cm,2.5 cm}]{geometry}
\usepackage{amsmath}
\usepackage{apacite}
\usepackage[square, numbers]{natbib}
\usepackage{amssymb}
\usepackage{graphicx}
\usepackage{algorithm}
\usepackage{algorithmic}
\usepackage{subfigure}
\usepackage{multicol}
\usepackage{multirow}
\usepackage{bbm}
\usepackage{url}
\usepackage{tikz,pgfplots}
\usepackage{amsthm}
\usepackage{authblk}
\usepackage[pdfborder={0 0 0}]{hyperref}
\newcommand{\Keywords}[1]{\par\noindent{\small{\em Keywords\/}: #1}}
\def\ZZ{\mathbb{Z}}

\def\RR{\mathbb{R}}

\def\EE{\mathbb{E}}
\def\PP{\mathbb{P}}

\def\T{{\cal T}}

\begin{document}
\onehalfspacing

\title{Imperfect demand estimation for new product production planning}

\author[1]{Antoine Deza}
\author[2]{Kai Huang}
\author[3]{Michael R. Metel}

\affil[1]{Advanced Optimization Laboratory, Department of Computing and Software, McMaster University, Hamilton, Ontario, Canada\\
\url{deza@mcmaster.ca}}
\affil[2]{DeGroote School of Business, McMaster University, Hamilton, Ontario, Canada\\
\url{khuang@mcmaster.ca}}
\affil[3]{School of Computational Science and Engineering, McMaster University, Hamilton, Ontario, Canada\\
\url{metelm@mcmaster.ca}}

\maketitle

\begin{abstract}
We are interested in the effect of consumer demand estimation error for new products in the context of production planning. An inventory model is proposed, whereby demand is influenced by price and advertising. The effect of parameter misspecification of the demand model is empirically examined in relation to profit and service level feasibility. Faced with an uncertain consumer reaction to price and advertising, we find that it is safer to overestimate rather than underestimate the effect of price on demand. Moreover, under a service level constraint it is safer to overestimate the effect of advertising, whereas for strict profit maximization, underestimating the effect of advertising is the conservative approach.
\end{abstract}


\section{Introduction}

In this paper we investigate the effect of model parameter estimation error to determine best practices when estimating the effect of price and advertising in relation to inventory management~\cite{zipkin}. We consider a single period inventory model with a minimum service level constraint~\cite{Chen2001}, with the objective of maximizing profit under consumer demand uncertainty.\\

The Bass Model~\cite{Bass69} is a differential equation which is widely used to forecast new product adoption. Extensions have been made, such as the Generalized Bass Model~\cite{bass1994}, which incorporates both price and advertising. For this paper we use an approximation of the  Piecewise-Diffusion Model (PDM) of Nui~\cite{Niu06}, which extends the original Bass Model by incorporating demand uncertainty, as well as price and advertising, resulting in a superior fit compared to the previous models in empirical testing.\\

In the paper of Lim et al.~\cite{lim2013}, the misestimation of supply chain disruption probabilities was investigated. It was found that overestimating disruption probabilities reduces the expected cost when compared to underestimation. When faced with estimation uncertainty, this presents the managerial insight that having a bias towards overestimation prevents excessive costs.
In the problem setting of this paper, the estimation uncertainty lies in the consumer demand model. It is unclear what the effect of the estimation error of the consumer demand's response to price and advertising is on profit and service level feasibility, and when faced with uncertainty, what the conservative approach to estimation would be. In an attempt to answer these questions, we conduct an empirical study, whereby the optimal solution is found for our inventory model under what is considered to be the true consumer demand dynamics, after which optimal solutions are found under biased responses to price and advertising to determine the effect of misestimation.\\

We present an overview of the PDM and briefly trace its roots in Section~\ref{sec:pdm}. In Section~\ref{sec:mod} the inventory optimization problem is presented as well as the formulation of its approximation which we solve. Included are details of the calibration of the model, the problem instances which we are interested in, as well as a justification of our approximation in terms of confidence intervals of the error. Details of the computational experiments are described in Section~\ref{sec:exp}, with a commentary on the results. The conclusions and future research directions are summarized in Section~\ref{sec:con}, with the results of the experiment graphically presented in the Appendix.

\section{Piecewise-Diffusion Model (PDM)}
\label{sec:pdm}

The Bass model proposes that the number of adopters through time, $N(t)$, can be modeled by the differential equation $\frac{d}{dt}N(t)=(m-N(t))(p+\frac{q}{m}N(t))$, where $m$ is the market size, $p$ is the coefficient of innovation, which is the consumer's  intrinsic desire to purchase the product, and $q$ is the coefficient of imitation, which models the influence of existing adopters on the consumer, whose solution is $N(t,m,p,q)=m\frac{1-e^{-(p+q)t}}{1+(q/p)e^{-(p+q)t}}$.\\

The Stochastic Bass Model (SBM) assumes that consumer adoption follows a pure birth process, where $A_m(t)$ is the cumulative number of adopters by time $t$. The transition rate from adoption $j$ to $j+1$ is $\lambda_{mj}=(m-j)(\alpha+\frac{\beta}{m-1}j)$, where $\alpha$ and $\beta$, the intrinsic adoption rate and the induction rate, can be interpreted in the same manner as $p$ and $q$ in the Bass model. Let this be referred to as an SBM with specification $\{m,\alpha,\beta\}$. A central limit theorem is derived in~\cite{Niu06}, where it is proved that as $m\to\infty$, $\frac{A_m(t)-N(t,m,\alpha,\beta)}{\sqrt{\psi(t,m,\alpha,\beta)}}$ converges in distribution to a standard normal random variable, where $\psi(t,m,\alpha,\beta)=m\frac{(1+\beta/\alpha)e^{-2(\alpha+\beta)t}}{[1+(\beta/\alpha)e^{-(\alpha+\beta)t}]^4}\{e^{(\alpha+\beta)t}-1+2(\frac{\beta}{\alpha})(\alpha+\beta)t+(\frac{\beta}{\alpha})^2(1-e^{-(\alpha+\beta)t})\}$, so that for $m$ sufficiently large, we can approximate $A_m(t)$ as normal with mean $N(t,m,\alpha,\beta)$ and variance $\psi(t,m,\alpha,\beta)$.\\

The Piecewise Stochastic Bass Model assumes a sequence of time intervals where adoption levels will be observed. Let $a$ be the total number of adopters up to the present time. The model assumes that of the total available potential adopters $m-a$, only $(m-a)\pi$ are true prospects, where $\pi$ is the participation fraction, and the remainder are dormant. We can simulate the demand up to time $t$ under this formulation as an SBM with specification $\{(m-a)\pi,\hat{\alpha},\hat{\beta}\}$, where $\hat{\alpha}=\alpha+\frac{\beta}{m-1}a$ and $\hat{\beta}=\frac{(m-a)\pi-1}{m-1}\beta$.\\

The PDM incorporates the central limit theorem result, as well as an additional variance component $\delta^2$ to capture exogenous disturbance and model misspecification. The demand over time $t$ is approximated as normal with mean $\mu=N(t,(m-a)\pi,\hat{\alpha},\hat{\beta})$ and variance $\sigma^2=\psi(t,(m-a)\pi,\hat{\alpha},\hat{\beta})+\delta^2t$.
Using the PDM, we are able to influence future demand by the choice of the product price $p$ and advertising spending $v$.
Their effect is modeled by setting $\pi=\pi_m\{1-[(1-\frac{\pi_p}{\pi_m})e^{-\gamma_p v}]^{(\frac{p}{p_{ref}})^{-\eta}}\}$ and replacing $\beta$ by $\beta[1+\gamma_b (v_0+v)]$, where $\pi_m$ is the maximum possible participation fraction, $\pi_p$ is the value of $\pi$ when $p=p_{ref}$, which is a calibration reference price, $\eta$ controls price sensitivity, $\gamma_p$ controls the impact of $v$, and $\gamma_b$ scales the increase in influence of existing adopters from the aid of the advertising over the product's life.\\

The sales trajectory of room air conditioners from 1949-1961, Table $\ref{T2}$ of the Appendix, was used in the empirical study in~\cite{Niu06}, showcasing the superior fit of the PDM compared to the Bass and Generalized Bass Models, with a reduction in the sum of squared errors of 94.3\% and 84.4\% respectively. This dataset has been used extensively in the past, including the papers describing these two past models. The PDM was fit to the historical data using maximum likelihood estimation, while the latter two were fit using nonlinear-least squares.\footnote{For new products with no sales history, this process is not possible, which in part motivated this research.} In this paper we utilize the actual history parameterization, fit to this dataset, which is in Table $\ref{T1}$ of the Appendix.

\section{Optimization Model}
\label{sec:mod}
We consider an inventory model with zero lead time, variable ordering cost $c$, and salvage price $s<c$. At the beginning of the time period, we set the price $p$ of our product, we determine the amount of advertising spending $v$, a product order $o$ is placed and received, and then the consumer demand $D$ is realized. We want to maximize profit subject to satisfying $D$ with probability $1-\theta$. Our sales over the period will be $\min\{o,D\}$, with our excess supply equal to $\max\{o-D,0\}$. The optimization problem is as follows.
\begin{alignat}{6}
\label{1}
\tag{1}
&\max&&\text{ }\EE (p\min\{o,D\}+s\max\{o-D,0\}-co-v)\nonumber\\
&\mbox{s.t. }&&\PP(o-D\geq 0)\geq 1-\theta\nonumber\\
&&& p,v,o\geq0\nonumber
\end{alignat}

\subsection{Program Formulation}

We use a sample average approximation (SAA) to approximate the objective function~\cite{birge1997}. We approximate the expected sales $\EE(\min\{o,D\})$ as $\{\frac{1}{N}\sum_{j=1}^{N}r_j:r_j\leq o, r_j\leq \mu+\sigma z_j\}$, where $z_j$ is a standard normal sample, and the expected excess supply $\EE(\max\{o-D,0\})$ as $\{o-\frac{1}{N}\sum_{j=1}^{N}r_j\}$. The chance constraint can be written as
$o\geq\mu + \sigma\Phi^{-1}(1-\theta)$~\cite{boyd}, where $\Phi^{-1}$ is the inverse cumulative distribution function of the standard normal distribution. The objective contains bilinear terms, but for a fixed value of $p$, the objective of $(\ref{2})$ becomes linear, so we only consider a finite number of values for $p$.
\begin{alignat}{6}
\label{2}
\tag{2}
&\max&&\text{ }\frac{(p-s)}{N}\sum_{j=1}^{N}r_j+(s-c)o-v\nonumber\\
&\mbox{s.t. }&&o\geq\mu + \sigma\Phi^{-1}(1-\theta)\nonumber\\
&&&r_j\leq o\hspace{10 pt}j=1,...,N\nonumber\\
&&&r_j\leq \mu+\sigma z_j\hspace{10 pt}j=1,...,N\nonumber\\
&&&o\geq 0\nonumber\\
&&&p_{min}\leq p\leq p_{max}\nonumber\\
&&&0\leq v \leq v_{max}\nonumber\\
&&&p\in \ZZ.\nonumber
\end{alignat}

Note that $\mu$ and $\sigma$ are non-convex functions of $p$ and $v$. We approximate $\mu$ and $\sigma$ as piecewise linear functions using the logarithmic disaggregated convex combination (DLog) model of Vielma et al.~\cite{Vielma2010}. A Delaunay triangulation is used to segment the price and advertising domain into a set of triangles $\T$. DLog requires $\lceil log_2|\T|\rceil$ binary variables, enforcing a convex combination of the vertices of a single triangle to represent the values of $\mu$ and $\sigma$.

\subsection{Inventory Scenarios and Model Calibration}

We considered product costs of $60\%$ and $80\%$ of the historical 1949 price of a room air conditioner, \$410, $c_1=246$ and $c_2=328$, with $p_{min}=350$, $p_{max}=450$, $v_{max}=100$, and $s=0.1c$. Our test instances consist of values of $\theta=\{1,0.5,0.25,0.05\}$.\\

Our piecewise linear approximation of $\mu$ and $\sigma$ was constructed in the following manner. Beginning with the extreme points of our domain as vertices, we iteratively added vertices by taking a Delaunay triangulation of the current vertex set and finding the triangle with the centroid with the largest Euclidian norm of the percentage error of $\mu$ and $\sigma$ and their piecewise approximations, $\hat{\mu}$ and $\hat{\sigma}$. This error was then compared to the error of the midpoint of each edge of the triangle, with the point with the largest error added to the vertex set, with $p$ rounded to the nearest integer. The approximation was limited to the use of 10 binary variables. Taking 20,000 samples to construct empirical distribution functions of the percentage error of $\hat{\mu}$ and $\hat{\sigma}$, confidence intervals were found using the Dvoretzky-Keifer-Wolfowitz inequality~\cite{massart1990}, which states that for an empirical distribution function with $n$ samples, $F_n(x)$, and for any $x\in \RR$, $\PP(F_n(x)-F(x)>\epsilon)\leq e^{-2n\epsilon^2}$ for every $\epsilon\geq\sqrt{\frac{1}{2n}\ln2}$. This implies that $F(x)\geq(F_n(x)-\epsilon)(1-e^{-2n\epsilon^2})$. Taking a value of $\epsilon=0.014$, confidence intervals were found to be $\PP(|\frac{\mu-\hat{\mu}}{\mu}|\leq 0.0024)\geq0.95$ and $\PP(|\frac{\sigma-\hat{\sigma}}{\sigma}|\leq 0.00050)\geq0.95$.\\

The SAA sample size, $N=15,000$, was chosen to ensure the sample problem optimal objective $z^*_N$ is close to the true optimal objective value $z^*$ with high probability. We consider the convergence of the most challenging problem, namely, when $c=c_2$ and $\theta=0.05$. We are interested in bounding the error due to sampling, so let $z(p,v,o)=(p-s)(\hat{\mu} F(o)-\hat{\sigma}^2f(o))+po(1-F(o))+soF(o)-co-v$, where $F(x)$ and $f(x)$ are the cumulative and probability distribution functions of $D\sim N(\hat{\mu},\hat{\sigma}^2)$, which is the objective value of $(\ref{1})$ using $\hat{\mu}$ and $\hat{\sigma}$. A confidence interval for the optimality gap was calculated based on the technique of Mak et al.~\cite{Mak1999}. $M=20$ instances of $(\ref{2})$ were solved with optimal values $z^{*i}_N$ and objective values of $z(p^*_i,v^*_i,o^*_i)$. Let $\hat{\mu}^N_Z$ and $\hat{\sigma}^N_Z$, and $\hat{\mu}_Z$ and $\hat{\sigma}_Z$, equal the sample mean and standard deviations of $z^{*i}_N$ and $z(p^*_i,v^*_i,o^*_i)$, respectively. When solving a single instance of $(\ref{2})$ the $1-\alpha$ confidence interval of the optimality gap is estimated as $[\hat{\mu}_Z+t_{\frac{\alpha}{2},M-1}\hat{\sigma}_Z\leq z^*\leq \hat{\mu}^N_Z + t_{1-\frac{\alpha}{2},M-1}\hat{\sigma}^N_Z]$, where $t_{\frac{\alpha}{2},M-1}$ is the $\frac{\alpha}{2}$-critical value of the t-distribution with $M-1$ degrees of freedom. The percentage error optimality gap confidence interval with $\alpha=0.05$ was found to be $\PP(\frac{z^{*}_N-z(p,v,o)}{z(p,v,o)}\leq 0.0183)\geq 0.95$. Virtually all of the error came from the $z^{*i}_N$, as each problem instance found the same optimal $p^*_i$, and the same $v^*_i$ and $o^*_i$ up to 15 significant digits.

\section{Computational Experiments}
\label{sec:exp}
We are interested in the effect of parameter misspecification on profit and feasibility. In particular, if the effect is asymmetrical, this gives guidance when having to estimate consumer behaviour for new products with no prior history. We observe the effect of underestimating and overestimating the influence of price, $\eta$, and the influence of advertising, $\gamma_b$ and $\gamma_p$, which we denote simply as $\gamma$. Given what we consider a true parameter value $x_0$ from Table \ref{T1}, we repeated the process described in the previous section, estimating $\mu$ and $\sigma$ by a piecewise linear function and solving (\ref{2}) for the optimal values $p^*$, $v^*$, and $o^*$ for $x=\{-0.6x_0,-0.3x_0,x_0,0.3x_0,0.6x_0\}$, then the expected profit and the feasibility assuming $x_0$ was observed. All computing was conducted on a Windows 7 Home Premium 64-bit, Intel Core i5-2320 3GHz processor with 8 GB of RAM. The implementation was done in Matlab R2012a interfaced with Gurobi 6.0 using YALMIP~\cite{lof04} dated November 27, 2014. The results are shown graphically in Figures 1 to 4. Figure 1 displays the profit when optimizing over different values of $\eta$. Figure 2 displays the optimal order quantity $o^*$ in relation to the minimum order quantity $m$ required for feasibility, presented as a percentage difference, $\frac{o^*-m}{m}\times 100$. Figures 3 and 4 display the same for varying $\gamma$.\\

When $\eta$ is underestimated, the price is increased to take advantage of the subdued decrease in demand. As a result, too much is ordered, significantly decreasing profit. When $\eta$ is overestimated, price is decreased slightly with an assumed exaggerated increase in demand, again resulting in an excessive order given the price. An interesting observation is regardless of over or underestimating $\eta$, the solution is feasible. So with the focus now only on profit, from Figure 1, we conclude that it is better to err on the side of overestimating the effect of price on consumer demand, which decreases profit at a lower rate than underestimating.\\

When $\gamma$ is underestimated, the effect of advertising on demand in underestimated, resulting in an insufficient quantity of product ordered and an infeasible solution. The opposite effect occurs when $\gamma$ is overestimated, resulting in an excessive, but feasible order. With no regard to service levels, we observe from Figure 3 that it is more profitable to underestimate rather than overestimate the value of $\gamma$, but given a service level constraint, underestimation causes infeasibility, whereas overestimation ensures feasibility.

\section{Conclusions}
\label{sec:con}

This paper has examined the effect of over and underestimating the influence of price and advertising on consumer demand in the context of production planning. This is of particular interest for a new product with no prior sales history to aid in decision making. From an empirical study, we have found that the error is asymmetrical. Faced with uncertainty, it is prudent to overestimate the effect of price, resulting in a lower rate of loss in profit. Underestimating the effect of advertising results in a superior profit, but given a service level constraint, it will result in an insufficient order quantity, whereas overestimating the effect of advertising will ensure feasibility. We see the potential for future work stemming from this paper. We have focused on a single period model in order to capture the relationships between demand factors and profit and feasibility as clearly as possible, but the extension to a multi-stage inventory model with service level constraints would be interesting from a modeling and computational aspect. We examined the two factors to consumer demand which we felt are of most interest to business managers, but perhaps future research could examine the effect of misestimating other factors, such as the maximum participation fraction $\pi_m$, which is closely related to the estimation of the market size.

\bibliographystyle{abbrvnat}
\bibliography{SBMInv}

\section*{Appendix}

\begin{table}[htb]
\centerline{
\resizebox{1.25\textwidth}{!}{
\begin{tabular}{lrrrrrrrrrrrrr}
	\hline\noalign{\vskip 2pt}
Year               & 1949 & 1950 & 1951 & 1952 & 1953 & 1954  & 1955 & 1956 & 1957  & 1958 & 1959 & 1960 & 1961 \\\noalign{\vskip 2pt}
 \hline\noalign{\vskip 2pt}
Sales (M)          & 96   & 195  & 238  & 365  & 1045 & 1230  & 1270 & 1828 & 1586  & 1673 & 1660 & 1580 & 1500 \\
Price (\$)         & 410  & 370  & 365  & 388  & 335  & 341   & 320  & 293  & 310   & 279  & 269  & 275  & 259  \\
Advertising (\$MM)  & 0    & 0.615& 1.198& 3.196& 5.34 & 14.372& 9.391& 13.61& 16.785& 9.238& 5.863& 3.923& 1.493\\\noalign{\vskip 2pt}
 \hline
\end{tabular}}}
\caption{Room air conditioner data from 1949-1961}\label{T2}
\end{table}

\begin{table}[htb]
\centerline{
\resizebox{0.9\textwidth}{!}{
\begin{tabular}{cccccccccc}
\hline\noalign{\vskip 2pt}    	
m ($10^3$) & $a_0$ ($10^3$) & $\pi_p$ & $\alpha$ & $\beta$ & $\delta$ & $\eta$ & $\pi_m$ & $\gamma_p$ & $\gamma_b$ \\\noalign{\vskip 2pt}    \hline\noalign{\vskip 2pt}
53,291  & 744  & 0.005191  & 0  & 19.14  & 39.52   & 6.218  & 0.04195  & 0.009746 & 0.3704\\\noalign{\vskip 2pt}
\hline
\end{tabular}}}
\caption{Actual history parameterization}\label{T1}
\end{table}

 \begin{figure}[htb]
 \centerline{
\resizebox{1.0\textwidth}{!}{
\begin{tikzpicture}
\scriptsize
    \begin{axis}[xlabel=Percent change of $\eta$,ylabel=Profit,
                 legend style={legend pos=outer north east}]
                \addplot[mark=*,mark size=0.75,mark options={color=orange},draw=orange]
                table[x=x,y=y600]
            {results8.dat};
                \addplot[mark=*,mark size=0.75,mark options={color=purple},draw=purple]
                table[x=x,y=y6050]
            {results8.dat};
                \addplot[mark=*,mark size=0.75,mark options={color=teal},draw=teal]
                table[x=x,y=y6075]
            {results8.dat};
                \addplot[mark=*,mark size=0.75,mark options={color=blue},draw=blue]
                table[x=x,y=y6095]
            {results8.dat};
                \addplot[mark=*,mark size=0.75,mark options={color=brown},draw=brown]
                table[x=x,y=y800]
            {results8.dat};
                \addplot[mark=*,mark size=0.75,mark options={color=magenta},draw=magenta]
                table[x=x,y=y8050]
            {results8.dat};
                \addplot[mark=*,mark size=0.75,mark options={color=black},draw=black]
                table[x=x,y=y8075]
            {results8.dat};
                \addplot[mark=*,mark size=0.75,mark options={color=cyan},draw=cyan]
                table[x=x,y=y8095]
            {results8.dat};
    \legend{$c=c_1$ $\theta=1\text{ }\text{ }\text{ }$,$c=c_1$ $\theta=0.5\text{ }$,$c=c_1$ $\theta=0.25$,$c=c_1$ $\theta=0.05$,$c=c_2$ $\theta=1\text{ }\text{ }\text{ }$,$c=c_2$ $\theta=0.5\text{ }$,$c=c_2$ $\theta=0.25$,$c=c_2$ $\theta=0.05$}
    \end{axis}
\end{tikzpicture}}}
\caption{Expected profit when optimizing over different values of $\eta$.} \label{T10}
\end{figure}

 \begin{figure}[htb]
 \centerline{
\resizebox{1.0\textwidth}{!}{
\begin{tikzpicture}
\scriptsize
    \begin{axis}[xlabel=Percent change of $\eta$,ylabel=\% from minimum feasible order,
                 legend style={legend pos=outer north east}]
                \addplot[mark=*,mark size=0.75,mark options={color=purple},draw=purple]
                table[x=x,y=y6050]
            {results16.dat};
                \addplot[mark=*,mark size=0.75,mark options={color=teal},draw=teal]
                table[x=x,y=y6075]
            {results16.dat};
                \addplot[mark=*,mark size=0.75,mark options={color=blue},draw=blue]
                table[x=x,y=y6095]
            {results16.dat};
                \addplot[mark=*,mark size=0.75,mark options={color=magenta},draw=magenta]
                table[x=x,y=y8050]
            {results16.dat};
                \addplot[mark=*,mark size=0.75,mark options={color=black},draw=black]
                table[x=x,y=y8075]
            {results16.dat};
                \addplot[mark=*,mark size=0.75,mark options={color=cyan},draw=cyan]
                table[x=x,y=y8095]
            {results16.dat};
                \addplot[draw=gray, dashed]
                table[x=x,y=zero]
            {results16.dat};
    \legend{$c=c_1$ $\theta=0.5\text{ }$,$c=c_1$ $\theta=0.25$,$c=c_1$ $\theta=0.05$,$c=c_2$ $\theta=0.5\text{ }$,$c=c_2$ $\theta=0.25$,$c=c_2$ $\theta=0.05$}
    \end{axis}
\end{tikzpicture}}}
\caption{Percentage from minimum feasible order when optimizing over different values of $\eta$.} \label{T11}
\end{figure}

 \begin{figure}[htb]
 \centerline{
\resizebox{1.0\textwidth}{!}{
\begin{tikzpicture}
\scriptsize
    \begin{axis}[xlabel=Percent change of $\gamma$,ylabel=Profit,
                legend style={legend pos=outer north east}]
                \addplot[mark=*,mark size=0.75,mark options={color=orange},draw=orange]
                table[x=x,y=y600]
            {results9.dat};
                \addplot[mark=*,mark size=0.75,mark options={color=purple},draw=purple]
                table[x=x,y=y6050]
            {results9.dat};
                \addplot[mark=*,mark size=0.75,mark options={color=teal},draw=teal]
                table[x=x,y=y6075]
            {results9.dat};
                \addplot[mark=*,mark size=0.75,mark options={color=blue},draw=blue]
                table[x=x,y=y6095]
            {results9.dat};
                \addplot[mark=*,mark size=0.75,mark options={color=brown},draw=brown]
                table[x=x,y=y800]
            {results9.dat};
                \addplot[mark=*,mark size=0.75,mark options={color=magenta},draw=magenta]
                table[x=x,y=y8050]
            {results9.dat};
                \addplot[mark=*,mark size=0.75,mark options={color=black},draw=black]
                table[x=x,y=y8075]
            {results9.dat};
                \addplot[mark=*,mark size=0.75,mark options={color=cyan},draw=cyan]
                table[x=x,y=y8095]
            {results9.dat};
    \legend{$c=c_1$ $\theta=1\text{ }\text{ }\text{ }$,$c=c_1$ $\theta=0.5\text{ }$,$c=c_1$ $\theta=0.25$,$c=c_1$ $\theta=0.05$,$c=c_2$ $\theta=1\text{ }\text{ }\text{ }$,$c=c_2$ $\theta=0.5\text{ }$,$c=c_2$ $\theta=0.25$,$c=c_2$ $\theta=0.05$}
    \end{axis}
\end{tikzpicture}}}
\caption{Expected profit when optimizing over different values of $\gamma$.} \label{T12}
\end{figure}

 \begin{figure}[htb]
 \centerline{
\resizebox{1.0\textwidth}{!}{
\begin{tikzpicture}
\scriptsize
    \begin{axis}[xlabel=Percent change of $\gamma$,ylabel=\% from minimum feasible order,legend style={legend pos=outer north east}]
                \addplot[mark=*,mark size=0.75,mark options={color=purple},draw=purple]
                table[x=x,y=y6050]
            {results15.dat};
                \addplot[mark=*,mark size=0.75,mark options={color=teal},draw=teal]
                table[x=x,y=y6075]
            {results15.dat};
                \addplot[mark=*,mark size=0.75,mark options={color=blue},draw=blue]
                table[x=x,y=y6095]
            {results15.dat};
                \addplot[mark=*,mark size=0.75,mark options={color=magenta},draw=magenta]
                table[x=x,y=y8050]
            {results15.dat};
                \addplot[mark=*,mark size=0.75,mark options={color=black},draw=black]
                table[x=x,y=y8075]
            {results15.dat};
                \addplot[mark=*,mark size=0.75,mark options={color=cyan},draw=cyan]
                table[x=x,y=y8095]
            {results15.dat};
    \legend{$c=c_1$ $\theta=0.5\text{ }$,$c=c_1$ $\theta=0.25$,$c=c_1$ $\theta=0.05$,$c=c_2$ $\theta=0.5\text{ }$,$c=c_2$ $\theta=0.25$,$c=c_2$ $\theta=0.05$}
    \end{axis}
\end{tikzpicture}}}
\caption{Percentage from minimum feasible order when optimizing over different values of $\gamma$.} \label{T13}
\end{figure}

\end{document}